\documentclass[oneside,english,reqno]{amsart}
\usepackage{palatino}
\usepackage[T1]{fontenc}
\usepackage[latin1]{inputenc}
\usepackage{amssymb}
\IfFileExists{url.sty}{\usepackage{url}}
                      {\newcommand{\url}{\texttt}}

\makeatletter

\providecommand{\LyX}{L\kern-.1667em\lower.25em\hbox{Y}\kern-.125emX\@}

 \theoremstyle{remark}    
 \newtheorem*{acknowledgement*}{Acknowledgement} 
 \theoremstyle{plain}    
 \newtheorem*{cor*}{Corollary}
 \theoremstyle{definition}
 \newtheorem*{defn*}{Definition}
 \theoremstyle{plain}    
 \newtheorem{lem}{Lemma} 
 \theoremstyle{remark}
 \newtheorem*{rem*}{Remark}
 \theoremstyle{plain}    
 \newtheorem{thm}{Theorem} 

\usepackage{babel}
\makeatother
\begin{document}

\title{Excited random walk in three dimensions has positive speed}

\author{Gady Kozma}

\maketitle

\section{introduction}

Excited random walk is a model of a random walk on $\mathbb{Z}^{d}$
which, whenever it encounters a new vertex it receives a push toward
a specific direction, call it the {}``right'', while when it reaches
a vertex it {}``already knows'', it performs a simple random walk.
This model has been suggested in \cite{BW} and had since got lots
of attention, see \cite{V,Z}. The reason for the interest is that
it is situated very naturally between two classical models: random
walk in random environment and reinforced random walk. A reinforced
random walk is a walk on a graph (say $\mathbb{Z}^{2}$) that, whenever
it passes through an edge, it changes the weight of this edge, usually
positively (i.e.~the edge has now a greater probability to be chosen
when the random walk rereaches one of its end points) but possibly
also negatively, with the extreme being the {}``bridge-burning random
walk'' that can never traverse the same edge twice. The problem appears
naturally in brain research in connection with the evolution of neural
networks. Reinforced random walk models are notoriously difficult
to analyze, and even the question whether the simplest one-reinforced
random walk on $\mathbb{Z}^{2}$ is recurrent or transient is open.
See \cite{K90,KR99,PV99,DKL02} for some known results.

Random walk in a random environment is also a model in which the environment
is random, but independently of the walk. For example, one may throw
a coin at every point of $\mathbb{Z}$ to decide if at this point
the walk will have a push to the left or to the right, and then perform
random walk on the resulting weighted graph. The independence of the
walk from the environment turns out to be a powerful leverage, and
many very precise results are known. See e.g.~the book \cite{H95}.

Excited random walk has, seemingly, all the difficulties of reinforced
random walk: the environment depends on the walk, and in a dynamic
way. However, it has two significant advantages. The first is the
inherent directedness: the drift of excited random walk is always
in the same direction, and in particular, it can be coupled with simple
random walk so that the excited is always to the right of the simple
random walk. The second is the projected simple random walk of lower
dimension. Thus, for example, for the excited random walk in three
dimensions, its projection on the two directions orthogonal to our
{}``right'' is a simple two-dimensional random walk, up to a time
change.

Thus, for example, it is clear that the three dimensional excited
random walk is transient. Indeed, since a two-dimensional simple random
walk visits an order of $n/\log n$ vertices, the three dimensional
excited random walk must visit at least $n/\log n$ vertices. This
means, roughly, that $R(n)_{1}>n/\log n-C\sqrt{n\log \log n}$ ($x_{1}$
denoting the first, {}``left-right'' coordinate of $x$), and in
particular that $R(n)$ drifts to the right and returns to every point
only a finite number of times (in the two dimensional case this argument
does not work --- see \cite{BW} for a proof of this fact). The purpose
of this note is to improve this obvious remark. We shall show that
the factor $1/\log n$ is only an artifact of this argument, namely
we shall prove

\begin{thm}
\label{thm:exc}Let $R(n)$ be an $\epsilon $-excited random walk.
Then\[
\liminf _{n\to \infty }\frac{R(n)_{1}}{n}>0.\]

\end{thm}
The corresponding problem in two dimensions remains open. 

I believe that the lower limit above is in effect a limit. We will
not prove it, but our techniques show great independence between different
parts of an excited random walk, therefore it stands to reason that
it shouldn't be difficult.

An important element of the proof is a two dimensional result which
might be of independent interest --- indeed we already have another
application for it, \cite{ABK}. It reads

\begin{thm}
\label{thm:expsqrtlog}Let $R_{1}$ and $R_{2}$ be two independent
simple random walks on $\mathbb{Z}^{2}$ starting from $0$, $R_{1}$
of length $n$ and $R_{2}$ of length $m=\exp (\log ^{\mu }n)$ for
some $\mu \in [\frac{1}{2},1]$. Then\begin{eqnarray}
\lefteqn{\mathbb{P}(\#\{v:\exists t\leq m,\, R_{2}(t)=v\textrm{ but }\forall s\leq n,\, R_{1}(s)\neq v\}\leq m^{3/4})\leq } &  & \nonumber \\
 & \qquad \qquad \qquad \qquad \qquad \qquad  & \leq C\exp (-c\log ^{2\mu -1}n).\label{eq:thm1}
\end{eqnarray}

\end{thm}
(note that $\mu $ has the elegant expression $\mu =\log \log m/\log \log n$).
In particular, it shows that $R_{1}$ has at least $m^{3/4}$ {}``holes''
in a $\sqrt{m}$ vicinity of zero --- if $R_{2}$ finds them then
they must exist! The theorem is sharp in the following sense: with
probability $>c\exp (-C\log ^{2\mu -1}n)$, $R_{1}$ covers all of
$B(0,m)$. The (easy) proof of sharpness will only be sketched below.

\subsection{Around the proof}

The proof of theorem \ref{thm:exc} is only three pages long, but
its inductive nature, the number of parameters and their interdependencies
make it somewhat opaque. Therefore I feel compelled to make some vague
comments in preparation for the actual argument. The basic argument
is a block decomposition. This approach has been tried before, but
a straightforward attack does not work. If you divide your time span
$[0,n]$ into blocks of length $k$ and allow yourself to {}``lose
a factor of $\frac{1}{\log k}$'' in each, you are left with the
following obstacle: once you have one really bad block (which \emph{will}
happen, if $n$ is large enough), you have difficulties to say anything
useful about the next block. And then about the block following it.
And so on. Hence $k$ cannot be independent of $n$ --- it has to
be \emph{at least} $\log n$ to get something. Thus the basic block
approach gives (roughly) $R(n)/n>c/\log \log n$, but not a constant. 

The argument here tries to work around this problem by a {}``restart
mechanism'', namely some way to continue after encountering a bad
block. This mechanism, roughly, throws away a big chunk in this case,
initializes the process from two dimensional considerations, and is
then forced to {}``pay'' just a little for the bad block, and of
course, they happen very rarely. The {}``big chunk'' above, denote
it by $l$, where $n\gg l\gg k$, is simply an intermediate size block.
Since multiple layers are needed to get an actual constant, the easiest
method to describe the structure is inductive. Thus the reader should
probably keep in mind, while reading the proof, that it really describes
a multi-layer structure where layer $i+1$ is used to restart the
estimates in the rare events that a block in the $i$th level failed.

\begin{acknowledgement*}
I wish to thank Itai Benjamini for many useful discussions.
\end{acknowledgement*}

\subsection{Excited random walk --- notations}

Let $0<\epsilon \leq \frac{1}{6}$. An $\epsilon $-excited random
walk (in three dimensions) is a random sequence $\left\{ R(n)\right\} _{n=0}^{\infty }$
of points in $\mathbb{Z}^{3}$ with the distribution defined as follows.
$R(0)=(0,0,0)$. Denote $R(n)=:(x_{1},x_{2},x_{3})$. Then

\begin{enumerate}
\item If $R(n)=R(m)$ for some $m<n$ ($R(n)$ is {}``visited''), then
$R(n+1)$ is one of the six neighbors of $R(n)$ in $\mathbb{Z}^{3}$
with probability $\frac{1}{6}$ each.
\item Otherwise ($R(n)$ is {}``new''), the probability is $\frac{1}{6}+\epsilon $
for $R(n+1)=(x_{1}+1,x_{2},x_{3})$ and $\frac{1}{6}-\epsilon $ for
$R(n+1)=(x_{1}-1,x_{2},x_{3})$. The other neighbors have probability
$\frac{1}{6}$ each.
\end{enumerate}
In both cases, the random choice is independent of the past, except
for the position $R(n)$ and whether the vertex is visited or new.

If $\mathfrak{V}\subset \mathbb{Z}^{3}$ is any set and $x\in \mathbb{Z}^{3}$
is a point, then an $\epsilon $-excited random walk starting from
$(x,\mathfrak{V})$ is an $\epsilon $-excited random walk such that
$R(0)=x$, and such that if $R(n)\in \mathfrak{V}$ then rule 1 above
is applies to it regardless of the past of $R$, i.e.~all vertices
in $\mathfrak{V}$ are considered {}``visited''.

\subsection{Standard notations}

The notations $C$ and $c$ relate to absolute constants, which may
be different from place to place. Sometimes we shall number them for
clarity. $C$ will usually pertain to constants which are {}``large
enough'' and $c$ to constants which are {}``small enough''. The
notation $x\approx y$ is a short hand for $cx\leq y\leq Cx$. The
notations $\ll $ and $\gg $ have no particular additional mathematical
content over $<$ and $>$. We only use them to stress that in a specific
point the estimate is very rough, and that's OK because it is enough
for our purposes.

For a subset $A\subset \mathbb{Z}^{d}$, we denote by $\partial A$
the \emph{inner} boundary, namely all vertices $v\in A$ with at least
one neighbor outside $A$.

For a number $x$, $\left\lfloor x\right\rfloor $ will denote the
largest integer $\le x$ and $\left\lceil x\right\rceil $ will denote
the smallest integer $\geq x$.

\section{Simple random walk in two dimensions --- the $e^{C\sqrt{\log n}}$
phenomenon.}

\begin{lem}
\label{lem:0_not_too_small}Let $k\in \mathbb{N}$ and let $r>e^{\sqrt{k}}$.
Let $x_{1},\dotsc ,x_{k}\in \partial B(0,r)$ and let $R_{1},\dotsc ,R_{k}$
be random walks starting from $x_{i}$ and stopped on $B(0,2r)$.
Then\[
\mathbb{P}(R_{i}(t)\neq 0\: \forall i=1,\dotsc ,k\forall t)>ce^{-C\sqrt{k}}.\]

\end{lem}
\begin{proof}
Let $a$ be the harmonic potential on $\mathbb{Z}^{2}$ (see e.g.~\cite{S76}).
Let $p_{i}$ be the probability that $R_{i}$ hits $0$ before exiting
$B(0,2r)$. Let $\tau $ be the stopping time when $R_{i}$ hits either
$0$ or $\partial B(0,2r)$. Since $a$ is harmonic outside $0$,
we get\[
a(x_{i})=\mathbb{E}R_{i}(\tau )=p_{i}a(0)+(1-p_{i})\mathbb{E}(R_{i}(\tau )\, |\, R_{i}(\tau )\in \partial B(0,er)).\]
Since $a(x_{i})=\frac{2}{\pi }\log r+O(1)$ and $a(x)=\frac{2}{\pi }\log (2r)+O(1)$
for any $x\in \partial B(0,er)$ (see~\cite[P12.3, page 124]{S76}),
and since $a(0)=0$ we get\[
p_{i}=\frac{\log 2}{\log 2r}+O\left(\log ^{-2}r\right)=\frac{\log 2}{\log r}+O(\log ^{-2}r).\]
Taking $k$'th power we get\begin{align*}
\mathbb{P}(R_{i}(t)\neq 0\: \forall i=1,\dotsc ,k;\, \forall t) & =\left(1-\frac{\log 2}{\log r}+O(\log ^{-2}r)\right)^{k}>\\
 & >\left(1-\frac{\log 2}{\sqrt{k}}+O(k^{-1})\right)^{k}>ce^{-(\log 2)\sqrt{k}}.\qedhere 
\end{align*}

\end{proof}
\begin{lem}
\label{lem:0_not_too_small_cond}Let $k\in \mathbb{N}$ and let $r>e^{C_{1}\sqrt{k}}$
for some $C_{1}$ sufficiently big. Let $x_{1}^{\textrm{in}},\dotsc ,x_{k}^{\textrm{in}}\in \partial B(0,r)$,
let $x_{1}^{\textrm{out}},\dotsc ,x_{k}^{\textrm{out}}\in \partial B(0,2r)$
and let $R_{1},\dotsc ,R_{k}$ be random walks starting from $x_{i}^{\textrm{in}}$
and conditioned to hit $B(0,2r)$ at $x_{i}^{\textrm{out}}$. Then\[
\mathbb{P}(R_{i}(t)\neq 0\: \forall i=1,\dotsc ,k;\, \forall t)>ce^{-\sqrt{k}}.\]

\end{lem}
\begin{proof}
This follows as lemma \ref{lem:0_not_too_small} when one remembers
the following fact: if $R$ is an (unconditioned) random walk starting
from $x^{\textrm{in}}\in \partial B(0,r)$ and stopped on $\partial B(0,2r)$,
and if $\mathbf{E}$ is any event that depends only on the portions
of $R$ inside $B(0,r)$, then\begin{equation}
\mathbb{P}(\mathbf{E})\approx \mathbb{P}(\mathbf{E}\, |\, R\textrm{ hits }B(0,2r)\textrm{ in some }x^{\textrm{out}}).\label{eq:condsowhat}\end{equation}
This is well known. see e.g.~\cite[lemma A.5]{BK} (the result there
is for dimension $\geq 3$ but the same proof holds for dimension
$2$ with minimal changes). The constant $C_{1}$ comes from the constants
implicit in the $\approx $ notation in (\ref{eq:condsowhat}).
\end{proof}
\begin{defn*}
Let $R$ be a random walk with some stopping time $\tau $ and let
$B(x,r)$ be some ball. Define stopping times by $\tau _{0}^{\textrm{out}}=0$
and \begin{align}
\tau _{j}^{\textrm{in}} & :=\min \{t\geq \tau _{j-1}^{\textrm{out}}:R_{i}(t)\in B(x,r)\}\nonumber \\
\tau _{j}^{\textrm{out}} & :=\min \{t>\tau _{j}^{\textrm{in}}:R_{i}(t)\in \partial B(x,2r)\}.\label{eq:deftauj}
\end{align}
Let $J:=\max \{j:\tau _{j}^{\textrm{in}}<\tau \}.$ We call $J$ \textbf{the
number of visits to} $B(x,r)$. In many cases we will have $k$ walks
$R_{i}$ with stopping times $\tau _{i}$. In this case we define
$\tau _{i,j}^{\textrm{in/out}}$ and $J_{i}$ in the same manner,
and call $\sum J_{i}$ \textbf{the total number of visits to} $B(x,r)$.
Note that it is possible for (some of) the random walk to start inside
the ball $B(x,r)$. In this case $\tau _{1}^{\textrm{in}}=0$ and
this is considered the first visit.
\end{defn*}
\begin{lem}
\label{lem:xrc}Let $k\in \mathbb{N}$ and let $r>e^{C_{2}\sqrt{k}}$
for some $C_{2}$ sufficiently big. Let $x_{1}^{\textrm{in}},\dotsc ,x_{k}^{\textrm{in}}\in B(0,r)$,
$x_{1}^{\textrm{out}},\dotsc ,x_{k}^{\textrm{out}}\in B(0,2r)$ and
let $R_{1},\dotsc ,R_{k}$ be random walks starting from $x_{i}^{\textrm{in}}$
and conditioned to hit $B(0,2r)$ in $x_{i}^{\textrm{out}}$. Then\begin{equation}
\mathbb{P}(\#\{x\in B(0,r):R_{i}(t)\neq x\: \forall i=1,\dotsc ,k;\, \forall t\}>r^{7/4})>c.\label{eq:xrc}\end{equation}

\end{lem}
\begin{proof}
Clearly, we may assume $r$ is sufficiently large (in the sense that
$r>C$) and pay only in the constant $c$ in (\ref{eq:xrc}); and
since the probability is decreasing in $k$ that means we may also
assume $k$ is sufficiently large. Let $s=e^{C_{3}\sqrt{k}}$ for
some $C_{3}<C_{2}$ that will be fixed, together with $C_{2}$, only
later (however, the implicit constants $r_{\min }$ and $k_{\min }$
are assumed to be fixed after $C_{2}$ and $C_{3}$ and may depend
on their values). We do need to remark in this stage that $C_{2}-C_{3}$
is also {}``sufficiently large'' i.e.~during the proof we will
only add restrictions that increase it. Let $y_{1},\dotsc ,y_{n}\in B(0,r)$
satisfy that $B(y_{m},2s)$ are disjoint and that $B(y_{m},2s)\subset B(0,r)\setminus \{x_{1}^{\textrm{in}},\dotsc ,x_{k}^{\textrm{in}}\}$.
Clearly, we may assume $n>c(r/s)^{2}$. Examine the $i$'th walk
(for a while, everything below will depend on this $i$ but we will
not repeat this fact every time). Define stopping times $\tau _{j}^{\textrm{in/out}}=\tau _{i,j}^{\textrm{in/out}}$
as follows: $\tau _{0}^{\textrm{out}}=0$ and \begin{align}
\tau _{j}^{\textrm{in}} & :=\min \{t>\tau _{j-1}^{\textrm{out}}:\exists m,\, R_{i}(t)\in B(x_{m},s)\}\nonumber \\
\tau _{j}^{\textrm{out}} & :=\min \{t>\tau _{j}^{\textrm{in}}:R_{i}(t)\in \partial B(x_{m},2s)\}.\label{eq:deftauij}
\end{align}
We have left out the $\exists m$ which might be formally needed in
the definition of $\tau _{j}^{\textrm{out}}$ since, clearly, the
same $m$ holds for both $\tau _{j}^{\textrm{in}}$ and $\tau _{j}^{\textrm{out}}$
--- $m$ may only change when $j$ changes. It is easy to see that
for some $c$ sufficiently small,\[
\mathbb{P}(\tau _{j}^{\textrm{out}}-\tau _{j}^{\textrm{in}}<cs^{2}\, |\, R[0,\tau _{j}^{\textrm{in}}])<\frac{1}{4}\]
(i.e.~the estimate holds independently of the past). Notice that
this uses (\ref{eq:condsowhat}) to overcome the conditioning over
the past. Hence, easily,\begin{equation}
\mathbb{P}(\tau _{j}^{\textrm{out}}<cjs^{2})<\frac{1}{4}\quad \forall j\label{eq:taujlarge}\end{equation}
(not necessarily the same $c$, of course). Let $\sigma =\sigma _{i}$
be the stopping time when $R_{i}$ exits $B(0,2r)$. It is well known
that $\sigma $ is approximately $r^{2}$, with an exponentially decreasing
tail, i.e. \begin{equation}
\mathbb{P}(\sigma >\lambda r^{2})<Ce^{-c\lambda }\label{eq:tausmall}\end{equation}
(the only difficulty is that $R$ is conditioned, and (\ref{eq:condsowhat})
doesn't apply. Again, we refer to a proof of a high dimensional analog
result, \cite[lemma A.8]{BK}). 

We now return the notation $i$. Define $J_{i}:=\max \{j:\tau _{i,j}^{\textrm{out}}<\sigma _{i}\}$.
The same argument that gave (\ref{eq:taujlarge}) will give, with
another sum,\begin{equation}
\mathbb{P}\Big (\sum _{i}\tau _{i,J_{i}}<cs^{2}\sum _{i}J_{i}\Big )<\frac{1}{4},\label{eq:sumtauijlarge}\end{equation}
while a sum over (\ref{eq:tausmall}) would give\begin{equation}
\mathbb{P}\Big (\sum _{i}\sigma _{i}>\lambda kr^{2}\Big )<Ce^{-c\lambda }.\label{eq:sumtausmall}\end{equation}
Picking $\lambda $ sufficiently large such that this would be $<\frac{1}{4}$,
and combining (\ref{eq:sumtauijlarge}) we get an estimate for $\sum J_{i}$:\begin{equation}
\mathbb{P}\left(\sum J_{i}\leq C_{4}k(r/s)^{2}\right)\geq \frac{1}{2}.\label{eq:Ebig}\end{equation}
Denote this event by $\mathbf{E}$. Let $\Xi $ denote the space of
vectors $(\gamma _{1},\dotsc ,\gamma _{k},y_{1,1}^{\textrm{in}},\dotsc ,\linebreak [0]y_{1,\gamma _{1}}^{\textrm{in}},y_{2,1}^{\textrm{in}},\dotsc ,y_{k,\gamma _{k}}^{\textrm{in}},y_{1,1}^{\textrm{out}},\dotsc ,y_{k,\gamma _{k}}^{\textrm{out}})$
where $\gamma _{1},\dotsc ,\gamma _{k}$ are integers and the $y_{i,j}^{\textrm{in/out}}$
are points in some $\partial B(x_{m},s)$ and $\partial B(x_{m},2s)$
respectively. For every $\xi \in \Xi $ denote by $\mathbf{E}_{\xi }$
the event that $J_{i}=\gamma _{i}$ and that $R_{i}(\tau _{i,j}^{\textrm{in/out}})=y_{i,j}^{\textrm{in/out}}$.
Since $\xi $ clearly determines whether $\mathbf{E}$ happened or
not, define $\Xi '\subset \Xi $ to be the collection of all $\xi $'s
such that $\mathbf{E}_{\xi }$ ensures $\mathbf{E}$. Let $\mathbf{X}_{m}$
be the event that $R_{i}(t)\neq x_{m}$ for all $i$ and all $t$
and let $\mathbf{X}:=\sum \mathbf{1}\{\mathbf{X}_{m}\}$.

If the total number of visits of some ball $B(x_{m},s)$ by the $R_{i}$'s
is $\leq Ck$, we may apply lemma \ref{lem:0_not_too_small_cond},
if only $s>e^{C_{1}\sqrt{Ck}}$, where $C_{1}$ is from lemma \ref{lem:0_not_too_small_cond},
i.e.~if $C_{3}>C\cdot C_{1}$, and we get \begin{equation}
\mathbb{P}(\mathbf{X}_{m}\, |\, \mathbf{E}_{\xi })>ce^{-\sqrt{k}}.\label{eq:XmExi}\end{equation}
Notice that $\sum J_{i}$ is, in effect, the sum over all balls $B(x_{m},s)$
of the total number of visits of the $R_{i}$'s. Since the number
of balls is $>c(r/s)^{2}$, and since $\mathbf{E}$ says that $\sum J_{i}\leq C_{4}k(r/s)^{2}$,
we get that for at least half of the balls the number of visits is
$\leq Ck$. We now choose $C_{3}$ so that (\ref{eq:XmExi}) is satisfied
for at least half of the balls. Further, conditioning by $\mathbf{E}_{\xi }$
all the balls $B(x_{m},2s)$ are independent and we get using standard
estimates for independent variables, for $k$ sufficiently large,\[
\mathbb{P}(\mathbf{X}<c(r/s)^{2}e^{-\sqrt{k}}\, |\, \mathbf{E}_{\xi })<\frac{1}{2}\quad \forall \xi \in \Xi '.\]
Together with (\ref{eq:Ebig}) we get\begin{align*}
\mathbb{P}(\mathbf{X}\geq c(r/s)^{2}e^{-\sqrt{k}}) & \geq \sum _{\xi \in \Xi '}\mathbb{P}(\mathbf{E}_{\xi })\mathbb{P}(\mathbf{X}\geq c(r/s)^{2}e^{-\sqrt{k}}\, |\, \mathbf{E}_{\xi })\geq \frac{1}{2}\sum _{\xi \in \Xi '}\mathbb{P}(\mathbf{E}_{\xi })=\\
 & =\frac{1}{2}\mathbb{P}(\mathbf{E})\geq \frac{1}{4}.
\end{align*}

Hence we only need to check when $c(r/s)^{2}e^{-\sqrt{k}}>r^{7/4}$,
but this happens, for $k$ sufficiently large, when $C_{2}>4+8C_{3}$,
so we may now choose $C_{2}$ and the lemma is proved. 
\end{proof}
\begin{rem*}
Clearly, the same proof yields the stronger estimate\[
\mathbb{P}(\#\{x\in B(0,r):R_{i}(t)\neq x\: \forall i=1,\dotsc ,k;\, \forall t\}>r^{2-\epsilon })>c(\epsilon ).\]
However, we will not need it here.
\end{rem*}
\begin{lem}
\label{lem:allballs}Let $k\in \mathbb{N}$, let $s\geq 1$ and let
$r=se^{\sqrt{k}}$. Let $x_{1}^{\textrm{in}},\dotsc ,x_{k}^{\textrm{in}}\in B(0,r)$
and $x_{1}^{\textrm{out}},\dotsc ,x_{k}^{\textrm{out}}\in \partial B(0,2r)$
and let $R_{1},\dotsc ,R_{k}$ be random walks starting from $x_{i}^{\textrm{in}}$
and conditioned to hit $B(0,2r)$ in $x_{i}^{\textrm{out}}$. Let
$y$ satisfy $B(y,2s)\subset B(0,r)$ and $|y-x_{i}^{\textrm{in}}|\geq \frac{1}{4}r$
for all $x_{i}^{\textrm{in}}$. Let $J$ denote the total number of
visits of the $R_{i}$'s to $B(y,s)$. Then\[
\mathbb{P}(J>\lambda k)\leq e^{-(c\lambda -C)\sqrt{k}}\quad \forall \lambda >0.\]

\end{lem}
\begin{proof}
The proof is a simple variation on the classic estimate for sums of
independent variables. Let $J_{i}$ be the number of visits of $R_{i}$
to $B(y,s)$. The same harmonic potential estimates as in lemma \ref{lem:0_not_too_small}
show that for any point in $\partial B(y,2s)$, the probability to
hit $\partial B(0,2r)$ before hitting $B(y,s)$ is $\geq c/\log (r/s)=c/\sqrt{k}$
(here we use (\ref{eq:condsowhat})). This means that each $J_{i}$
has an exponential distribution with the tail decreasing faster than
$(1-c/\sqrt{k})^{n}$, or in other words,\begin{equation}
\mathbb{E}(e^{c_{1}J_{i}/\sqrt{k}}\, |\, J_{i}>0)\leq C\label{eq:expJsqrtkot0}\end{equation}
for some $c_{1}$ sufficiently small. 

We now use the condition $|y-x_{i}^{\textrm{in}}|>\frac{1}{4}r$:
a second application of the harmonic potential argument shows that
$\mathbb{P}(J_{i}>0)\leq C/\sqrt{k}$ (again using (\ref{eq:condsowhat})).
Plugging this into (\ref{eq:expJsqrtkot0}) gives\[
\mathbb{E}(e^{c_{1}J_{i}/\sqrt{k}})\leq 1+C/\sqrt{k}\]
and since the various $J_{i}$'s are independent we get\[
\mathbb{E}(e^{c_{1}J/\sqrt{k}})=\prod _{i=1}^{k}\mathbb{E}(e^{c_{1}J_{i}/\sqrt{k}})\leq (1+C/\sqrt{k})^{k}<e^{C\sqrt{k}}.\]
Therefore, by Chebyshev's inequality,\[
\mathbb{P}(J>\lambda k)=\mathbb{P}(e^{c_{1}J/\sqrt{k}}>e^{c_{1}\lambda \sqrt{k}})\leq \frac{\mathbb{E}e^{c_{1}J/\sqrt{k}}}{e^{c_{1}\lambda \sqrt{k}}}<e^{(C-c_{1}\lambda )\sqrt{k}}.\qedhere \]

\end{proof}
\begin{rem*}
The value $\frac{1}{4}$ is of course arbitrary --- it can be replaced
by any $\mu >0$ but the constants $C$ and $c$ from the formulation
of the lemma depend on this $\mu $.
\end{rem*}
\begin{proof}
[Proof of theorem \ref{thm:expsqrtlog}]As usual, we assume $n$ is
sufficiently large, as we may. The constant $n_{\min }$ will be fixed
last, at the very end of the proof. In particular we assume $n>1$
so that we have no problem dividing with $\log n$. Let $r=\sqrt{m}/\log n$.
Clearly, \begin{equation}
\mathbb{P}(R_{2}[0,m]\subset B(0,r))\leq Ce^{-cm/r^{2}}=Ce^{-c\log ^{2}n}\ll C\exp (-c\log ^{2\mu -1}n).\label{eq:B1}\end{equation}
Denote this {}``bad'' event by $\mathbf{B}_{1}$. 

Examine the number of visits of $R_{1}$ to $B(0,r)$. Let $x\in \partial B(0,2r)$
be some point and let $S$ be a random walk starting from $x$, and
let $\tau $ be the first time when $S\in B(0,r)\cup \partial B(0,2n)$.
Clearly, if $S(\tau )\in \partial B(0,2n)$ then $\tau \geq 2n-2r>n$
and then the usual harmonic potential argument gives\[
\mathbb{P}(S[0,n]\cap B(0,r)=\emptyset )\geq \mathbb{P}(S(\tau )\in \partial B(0,n))\geq c/\log (2n/r)>c/\log n.\]
Clearly, this implies that if $R_{1}(t)\in \partial B(0,2r)$ then
with probability $>c/\log n$ this is the last visit of $R_{1}$ to
$B(0,r)$. Hence we see that the number of visits $J$ has an exponentially
decreasing tail, and in particular, for any constant $c_{2}>0$,\begin{equation}
\mathbb{P}(J>c_{2}\log ^{2\mu }n)\leq C\exp (-c\log ^{2\mu }n/\log n)=C\exp (-c\log ^{2\mu -1}n)\label{eq:B2}\end{equation}
where the various $c$'s depend on $c_{2}$. We shall fix $c_{2}$
later on. This bad event (denote it by $\mathbf{B}_{2}$) is the one
with the largest probability, and the reason that the factor $\exp (-C\log ^{2\mu -1}n)$
appears in (\ref{eq:thm1}). Let \begin{equation}
k:=\left\lfloor c_{2}\log ^{2\mu }n\right\rfloor .\label{eq:defkmu}\end{equation}
Define the stopping times $\tau _{j}^{\textrm{in/out}}$ using (\ref{eq:deftauj})
for the ball $B(0,r)$, and from now on we shall examine $R_{1}[0,\tau _{k}^{\textrm{out}}]$
instead of $R_{1}[0,n]$. Similarly define $\sigma $ to be the stopping
time when $R_{2}$ exits $B(0,r)$ and replace $R_{2}[0,m]$ with
$R_{2}[0,\sigma ]$. More precisely, we shall show that \begin{equation}
\mathbb{P}(\#\{v:\exists t\leq \sigma ,\, R_{2}(t)=v;\, \forall s\leq \tau _{k}^{\textrm{out}},\, R_{1}(s)\neq v\}\leq m^{3/4})\leq C\exp (-c\log ^{2\mu -1}n),\label{eq:thmtauk}\end{equation}
which will finish the theorem with (\ref{eq:B1}) and (\ref{eq:B2}).

Let therefore $\mathbf{E}_{\xi }$, $\xi =(x_{1}^{\textrm{in}},\dotsc ,x_{k}^{\textrm{in}},x_{1}^{\textrm{out}},\dotsc ,x_{k}^{\textrm{out}})$,
be the event that $R_{1}(\tau _{j}^{\textrm{in/out}})=x_{j}^{\textrm{in/out}}$
(notice that $x_{1}^{\textrm{in}}=0$). Conditioning by $\mathbf{E}_{\xi }$
we get $k$ independent walks, each one conditioned to exit $B(0,2r)$
at a given point. We wish to use lemma \ref{lem:allballs} with $s:=re^{-\sqrt{k}}$.
First we note that $r=\sqrt{m}/\log n>c(\epsilon )\exp ((\frac{1}{2}-\epsilon )\log ^{\mu }n)$
for any $\epsilon >0$ so \begin{equation}
s>c(\epsilon )\exp \left(\left(\frac{1}{2}-\sqrt{c_{2}}-\epsilon \right)\log ^{\mu }n\right)\label{eq:sbigmu}\end{equation}
 so for $c_{2}$ sufficiently small and $n$ sufficiently large we
get $s\geq 1$, and we may apply lemma \ref{lem:allballs} in a meaningful
way and get, for every $B(y,s)$ satisfying $d(y,x_{i}^{\textrm{in}})\geq \frac{1}{4}r$,\begin{equation}
\mathbb{P}(J_{y}>\lambda k\, |\, \mathbf{E}_{\xi })\leq e^{-(c\lambda -C)\sqrt{k}}\quad \forall \lambda >0\label{eq:onenJy}\end{equation}
where $J_{y}$ is the number of visits to $y$. Pick $\lambda $ sufficiently
large such that the probability above is $\leq e^{-3\sqrt{k}}$. 

Examine the set \[
\mathcal{Y}:=\left(B(0,{\textstyle \frac{3}{4}}r)\setminus B(0,{\textstyle \frac{1}{4}}r)\right)\cap \left\lceil 4s\right\rceil \mathbb{Z}^{2}.\]
The proof of the theorem will follow from the interactions of $R_{1}$
and $R_{2}$ with the balls $B(y,2s)$, $y\in \mathcal{Y}$ (note
that they are disjoint). As in the proof of lemma \ref{lem:xrc},
we want to be able to consider the events inside each $B(y,s)$ as
independent. Define therefore stopping times $\rho _{i,j}^{\textrm{in/out}}$
for $i=1,2$, similarly to (\ref{eq:deftauij}), i.e \begin{align*}
\rho _{i,j}^{\textrm{in}} & :=\min \{t\geq \rho _{i,j-1}^{\textrm{out}}:\exists y\in \mathcal{Y}\textrm{ s.t. }R_{i}(t)\in B(y,s)\}\\
\rho _{i,j}^{\textrm{out}} & :=\min \{t>\rho _{i,j}^{\textrm{in}}:R_{i}(t)\in \partial B(y,2s)\}.
\end{align*}
and define $\mathbf{F}_{\zeta }$, where $\zeta =(J_{1},J_{2},z_{1,1}^{\textrm{in}},\dotsc ,z_{2,J_{2}}^{\textrm{in}},z_{1,1}^{\textrm{out}},\dotsc ,z_{2,J_{2}}^{\textrm{out}})$,
to be the event that $R_{i}(\sigma _{i,j}^{\textrm{in/out}})=z_{i,j}^{\textrm{in/out}}$
and that $\rho _{1,J_{1}+1}^{\textrm{in}}>\tau _{k}$ and $\rho _{2,J_{2}+1}^{\textrm{in}}>\sigma $,
or in other words, that the total number of visits of $R_{i}$ to
the balls $B(y,s)$, $y\in \mathcal{Y}$, is $J_{i}$. Denote the
collection of these $\zeta $'s by Z. Examine first $R_{1}$.

\medskip{}
\noindent \textbf{1. $R_{1}$ and the balls $B(y_{i},2s)$.}
\smallskip{}

\noindent Since $x_{1}^{\textrm{in}}=0$ and all other $x_{i}^{\textrm{in}}\in \partial B(0,r)$
we get that $|y-x_{i}|\geq \frac{1}{4}r$ for all $i\leq k$ and $y\in \mathcal{Y}$.
Hence, since $\#\mathcal{Y}<C(r/s)^{2}=Ce^{2\sqrt{k}}$ with (\ref{eq:onenJy})
and our choice of $\lambda $ we get\begin{align}
\mathbb{P}(\exists y\in \mathcal{Y}\textrm{ s.t. }J_{y}>\lambda k\, |\, \mathbf{E}_{\xi }) & \leq Le^{-3\sqrt{k}}<Ce^{-\sqrt{k}}\leq C\exp (-c\log ^{\mu }n)\ll \nonumber \\
 & \ll C\exp (-c\log ^{2\mu -1}n)\label{eq:defB3}
\end{align}
Denote this event by $\mathbf{B}_{3}$. Let $Z'\subset Z$ be the
subset of all $\zeta $'s such that $\mathbf{F}_{\zeta }$ implies
$\neg \mathbf{B}_{3}$ (clearly, if $\mathbf{F}_{\zeta }$ happened
then we can calculate the number of visits to every $B(y,s)$ and
know whether $\mathbf{B}_{3}$ happened or not). Conditioning by $\mathbf{F}_{\zeta }$,
$\zeta \in Z'$ we get that all balls $B(y,s)$ are independent, and
we may use lemma \ref{lem:xrc} for every $B(y,s)$, if only $s>e^{C_{2}\sqrt{\lambda k}}$.
Remembering (\ref{eq:defkmu}) and (\ref{eq:sbigmu}), and comparing
the exponents, we see that this will hold (for $n$ sufficiently large)
if only \[
\frac{1}{2}-\sqrt{c_{2}}-\epsilon >C_{2}\sqrt{\lambda c_{2}}\]
which, again, holds if only $c_{2}$ is sufficiently small. The conclusion
of lemma \ref{lem:xrc} now reads

\begin{align}
 & \mathbb{P}(\#H_{y}>s^{7/4}\, |\, \mathbf{F}_{\zeta })>c\quad \forall y\in \mathcal{Y},\, \forall \zeta \in Z'.\label{eq:Bysepty}\\
 & H_{y}:=\{v\in B(y,s):\forall t<\tau _{k}^{\textrm{out}},\, R_{1}(t)\neq v\}.\nonumber 
\end{align}

\medskip{}
\noindent \textbf{2. $R_{2}$ and the balls $B(y_{i},2s)$.}
\smallskip{}

\noindent The last conclusion, (\ref{eq:Bysepty}), says in effect
that many balls $B(y,s)$ have {}``large $R_{1}$ holes'' (the $H_{y}$'s)
in them. Here we shall complement this with proving that $R_{2}$
passes thr\-ough many $B(y,s)$'s and at least in one of them, through
a sizable part of the $R_{1}$ hole. 

Easily, if $B(z,4s)\subset B(0,\frac{3}{4}r)\setminus B(0,\frac{1}{4}r)$,
and if $S$ if a random walk starting from $z$ and stopped at $\partial B(z,4s)$,
then \[
\mathbb{P}\left(S\cap \bigcup _{y\in \mathcal{Y}}B(y,s)=\emptyset \right)<1-c.\]
Hence the probability not to intersect $\bigcup B(y,s)$ has an exponentially
decreasing tail, as long as we are still within the annulus. In particular,
if $B(z,C_{5}\log n)\subset B(0,\frac{3}{4}r)\setminus B(0,\frac{1}{4}r)$
for some $C_{5}$ sufficiently large, and $S$ is stopped at $\partial B(z,C_{5}s\log n)$,
then\[
\mathbb{P}\left(S\cap \bigcup _{y\in \mathcal{Y}}B(y,s)=\emptyset \right)<\frac{1}{n^{2}}\]
Examine now an annulus $A(a):=B(0,a+(2C_{5}\log n+4)s)\setminus B(0,a)$
where $A(a)\subset B(0,\frac{3}{4}r)\setminus B(0,\frac{1}{4}r)$.
We get that with probability $>1-\frac{1}{n^{2}}$, $R_{2}$ intersects
a ball $B(y,s)\subset A(a)$. Taking $a_{i}=\frac{1}{4}r+i(2C_{5}\log n+4)s$
we get a sequence of $\approx r/s\log n$ disjoint annuli, and then
\[
\mathbb{P}(\forall a_{i}\, \exists y\textrm{ s.t. }B(y,s)\subset A(a)\textrm{ and }R_{2}\cap B(y,s)\neq \emptyset )>1-\frac{r}{ns\log n}\]
and in particular, if $\mathcal{Y}^{*}$ is the set of $y$'s such
that $R_{2}$ intersects $B(y,s)$, then\begin{equation}
\mathbb{P}(\#\mathcal{Y}^{*}<cr/s\log n)\ll C\exp (-c\log ^{2\mu -1}n).\label{eq:defB4}\end{equation}
Denote this event by $\mathbf{B}_{4}$.%
\footnote{This estimate is actually quite bad. The true expected value of $\#\mathcal{Y}^{*}$
is $(r/s)^{2}/\log (r/s)$, analogous to the fact that a random walk
of length $n$ passes through approximately $n/\log n$ distinct points.
However, it will do for our needs.%
} Define $Z''\subset Z'$ to be the subset of $\zeta $'s that ensure
that $\mathbf{B}_{4}$ did not happen.. 

\medskip{}
\noindent \textbf{3. The interaction between $R_{1}$ and $R_{2}$.}
\smallskip{}

Next examine one $B(y,s)$, $y\in \mathcal{Y}^{*}$. For every $v\in H_{y}$,
the harmonic potential argument shows that \[
\mathbb{P}(\exists t\leq \sigma :R_{2}(t)=v)\geq c/\log s\]
hence, if $h_{y}:=\#\{v\in H_{y}:\exists t\leq \sigma ,\, R_{2}(t)=v\}$
then $\mathbb{E}h_{y}\geq c\#H_{y}/\log s$, and of course $h_{y}\leq \#H_{y}$.
This shows that \[
\mathbb{P}(h_{y}>{\textstyle \frac{1}{2}}\mathbb{E}h_{y}\, |\, \mathbf{F}_{\zeta })>\frac{c}{\log s}\quad \forall y\in \mathcal{Y}^{*},\, \forall \zeta \in Z.\]
Remembering (\ref{eq:Bysepty}) and the independence of $R_{1}$ and
$R_{2}$ gives\[
\mathbb{P}(h_{y}>c_{3}s^{7/4}/\log s\, |\, \mathbf{F}_{\zeta })>\frac{c}{\log s}\quad \forall y\in \mathcal{Y}^{*},\, \forall \zeta \in Z'\]
for some $c_{3}$ sufficiently small. Remembering the definition
of $Z''$ (see below (\ref{eq:defB4})) we get\begin{align*}
\mathbb{P}(\forall y & \in \mathcal{Y}^{*},\, h_{y}\leq c_{3}s^{7/4}/\log s\, |\, \mathbf{F}_{\zeta })\leq \left(1-\frac{c}{\log s}\right)^{cr/s\log n}\\
 & \leq \exp \left(-ce^{\sqrt{k}}/\log ^{2}n\right)\ll C\exp (-c\log ^{2\mu -1}n)\quad \forall \zeta \in Z''
\end{align*}
(remember the definition of $k$, (\ref{eq:defkmu})). Throwing in
(\ref{eq:defB3}) and (\ref{eq:defB4}) and summing over all $\zeta $
we get\[
\mathbb{P}\Big (\sum _{y\in \mathcal{Y}^{*}}h_{y}\leq \frac{c_{3}s^{7/4}}{\log s}\Big )\leq C\exp (-c\log ^{2\mu -1}n).\]
However, this event is what we need in (\ref{eq:thmtauk})! Indeed,
directly from the definitions, \[
\#\{v:\exists t\leq \sigma ,\, R_{2}(t)=v;\, \forall s\leq \tau _{k}^{\textrm{out}},\, R_{1}(s)\neq v\}\geq \sum _{y\in \mathcal{Y}^{*}}h_{y},\]
so we need only explain why $m^{3/4}\leq c_{2}s^{7/4}/\log s$. Using
(\ref{eq:sbigmu}) we see that this is equivalent to, for $n$ sufficiently
large,\begin{equation}
\frac{7}{4}\left(\frac{1}{2}-\sqrt{c_{2}}-\epsilon \right)>\frac{3}{4},\label{eq:thhreefourth}\end{equation}
which holds for $c_{2}$ sufficiently small. Finally we may fix the
value of $c_{2}$, get (\ref{eq:thmtauk}) and hence the theorem.
\end{proof}
\theoremstyle{remark}

\newtheorem*{remarks*}{Remarks}

\begin{remarks*}

\begin{enumerate}
\item The only place in the proof the value $\frac{3}{4}$ appears is in
(\ref{eq:thhreefourth}). Hence, as in the remark following lemma
\ref{lem:xrc} (and using that remark), the theorem may be strengthened
to say\begin{eqnarray*}
\lefteqn{\mathbb{P}(\#\{v:\exists t\leq m,\, R_{2}(t)=v\textrm{ but }\forall s\leq n,\, R_{1}(s)\neq v\}\leq m^{1-\epsilon })\leq } &  & \\
 & \qquad \qquad \qquad \qquad \qquad \qquad  & \leq C(\epsilon )\exp (-c(\epsilon )\log ^{2\mu -1}n).
\end{eqnarray*}

\item As explained in the introduction, the theorem is sharp in the sense
that with probability $>c\exp (-C\log ^{2\mu -1}n)$, $R_{1}$ covers
all of $B(0,m)$. Roughly, the proof is as follows: lemma \ref{lem:0_not_too_small}
can be reversed to show that the probability to cover any point in
$B(0,\frac{1}{2}r)$ is $>1-C\exp (-c\sqrt{k})$. Similarly, the argument
leading to (\ref{eq:B2}) can be reversed to show that the probability
to have $\log ^{2\mu }n$ visits to $B(0,m)$ is $>c\exp (-C\log ^{2\mu -1}n)$,
and these two together give the result.
\end{enumerate}
\end{remarks*}

\section{Excited random walk in three dimensions}

The theorem will follow very easily from the following lemma. In effect,
the lemma is stronger than the theorem. The reason we need this stronger
formulation is its inductive proof.

\begin{lem}
\label{lem:excmain}Let $n\in \mathbb{N}$. Let $\mathfrak{V}\subset \left]-\infty ,-\left\lfloor n^{5/8}\right\rfloor \right]\times \mathbb{Z}^{2}$
be any configuration of visited vertices. Let $x\in \mathbb{Z}^{3}$.
Let $R$ be an $\epsilon $-excited random walk starting from $(x,\mathfrak{V})$
of length $2n$. Then\[
\mathbb{P}(R(n)_{1}>0\textrm{ and }R(2n)_{1}<R(n)_{1}+\alpha _{n}(\epsilon )n)\leq \exp (-c_{4}(\epsilon )\sqrt{\log n})\]
where the numbers $\alpha _{n}(\epsilon )$ satisfy a recursive condition
ensuring that $\alpha _{n}(\epsilon )\geq c(\epsilon )>0$.
\end{lem}
The mystery number $\frac{5}{8}$ is simply in the middle between
the $\frac{3}{4}$ of theorem \ref{thm:expsqrtlog} and $\frac{1}{2}$.
Since the $\frac{3}{4}$ of theorem \ref{thm:expsqrtlog} was an arbitrary
number $<1$, so is this $\frac{5}{8}$. For the impatient, the recursive
condition on the $\alpha (n)$ is (\ref{eq:defalphan}) below where
$k$ is defined in (\ref{eq:defknthm}) and where $\lambda $ is some
constant. It clearly ensures $\alpha _{n}\geq c$.

\begin{proof}
All the constants during the proof will depend on $\epsilon $, but
we will not repeat this fact and only write $C$ or $c$ instead of
$C(\epsilon )$ and $c(\epsilon )$. The lemma will be proved by induction,
so assume the lemma holds for any $k<n$ (we shall explain how to
deal with the case $n=1$, indeed with all sufficiently small $n$,
at the end). Due to this fact we need to pay special attention to
the constant $c_{4}$, to ensure that it is indeed a constant and
does not increase with $n$ --- hence none of the $C$ and $c$ below
will depend implicitly on $c_{4}$.

Our first observation is that one can couple (meaning, realizing them
on the same probability space) in the obvious way the excited random
walk in the interval $\left]n,2n\right]$ to a regular three dimensional
random walk $R'$ such that $R'(i)_{1}\leq R(n+i)_{1}$ for $i\leq n$.
For $R'$ we can use a simple estimate of binomial variables to say
that \begin{equation}
\mathbb{P}(\exists i\leq n:R'(i)_{1}<-n^{5/8})\leq Ce^{-cn^{1/8}}\ll C\exp (-c\sqrt{\log n}).\label{eq:Cirrelev}\end{equation}
Hence the same holds for $R\left]n,2n\right]$ and we conclude that
we do not need to know anything about $\mathfrak{V}$ --- with very
large probability, $R\left]n,2n\right]$ does not intersect $\mathfrak{V}$
and we need to investigate only its intersections with $R[0,n]$.

Let \begin{equation}
k:=\left\lceil n/\exp (\log ^{1/4}n)\right\rceil \label{eq:defknthm}\end{equation}
 and let $I_{i}:=\left]n+ik,n+(i+1)k\right]$ for $i=0,\dotsc ,\left\lceil \frac{n}{k}\right\rceil -1$.
Let \[
V_{i}:=\{v\in \mathbb{Z}^{3}:\exists t\in I_{i},\, R(t)=v\textrm{ and }\forall t\leq n+ik,\, R(t)\neq v\}.\]
Theorem \ref{thm:expsqrtlog} allows to estimate $\#V_{i}$ since
if the projections of the $V_{i}$'s on the second and third coordinates
are large then they themselves definitely will be. Denote the projection
by $P$. $P(R[0,n+ik]$) is a two dimensional random walk with the
length $\leq n+ik<2n$. $P(R(I_{i}))$ is a two dimensional random
walk whose length $m$ is a $(k,\frac{2}{3})$-binomial variable,
and in particular\[
\mathbb{P}(m<{\textstyle \frac{1}{2}}k)\leq Ce^{-ck}.\]
Assuming $m\geq \frac{1}{2}k$ we get $\mu =\log \log m/\log \log n=1-o(1)$
and hence theorem \ref{thm:expsqrtlog} says that, for $n$ sufficiently
large,\[
\mathbb{P}(\#V_{i}\leq m^{3/4})\leq C\exp (-c\log ^{2\mu -1}n)\leq C\exp (-c\sqrt{\log n})\quad \forall i.\]
 Examine the horizontal movement of $R$ during $I_{i}$ for one
$i$. There are $\leq k$ balanced horizontal moves (meaning that
they start from a visited vertex) and the number of unbalanced horizontal
moves is a $(\#V_{i},\frac{1}{3})$-binomial distribution. Hence with
probability $>1-Ce^{-ck}$ it is $\geq \frac{1}{6}\#V_{i}$ and if
$\#V_{i}>\frac{1}{2}k^{3/4}$ then we get a positive drift, namely
\begin{equation}
\mathbb{P}(R(n+(i+1)k)_{1}-R(n+ik)_{1}<c\epsilon k^{3/4})\leq C\exp (-c\sqrt{\log n})\quad \forall i\label{eq:RIilarge}\end{equation}
(in this formulation we no longer need to assume that $n$ is sufficiently
large).

Our purpose is to use the lemma inductively for every $I_{i}$, $i>0$.
Let therefore \[
a_{i}:=\max _{t\leq n+(i-1)k}R(t)_{1}+\left\lfloor k^{5/8}\right\rfloor .\]
We translate by $-a_{i}$ and use the lemma for $k$, which will now
read as\begin{equation}
\mathbb{P}(R(n+ik)_{1}>a_{i}\textrm{ and }R(n+(i+1)k)_{1}<R(n+ik)_{1}+\alpha _{k}k)\leq \exp (-c_{4}\sqrt{\log k}).\label{eq:defBi}\end{equation}
Denote this event by $\mathbf{B}_{i}$. The estimate (\ref{eq:defBi})
holds for any value of $R[0,n+(i-1)k]$ and hence we may rewrite it
as\[
\mathbb{P}(\mathbf{B}_{i}\, |\, \mathbf{B}_{0},\dotsc ,\mathbf{B}_{i-2})\leq \exp (-c_{4}\sqrt{\log k}).\]
Hence the sequence $\mathbf{B}_{2i}$ dominates a sequence of random
independent variables. Let $\lambda =\lambda (\epsilon )$ be an integer
parameter which we shall fix later. The simplest estimate now gives\[
\mathbb{P}\Big (\#\{i:B_{2i}\}\geq \lambda )\Big )\leq \left\lceil \frac{n}{k}\right\rceil ^{\lambda }e^{-\lambda c_{4}\sqrt{\log k}}\leq C\exp (-c_{4}\lambda (1-o(1))\sqrt{\log n}).\]
A similar calculation holds for the odd $\mathbf{B}_{2i+1}$, and
we get\begin{equation}
\mathbb{P}(\#\{i:\mathbf{B}_{i}\}\geq 2\lambda +2)\leq C\exp (-c_{4}\lambda (1-o(1))\sqrt{\log n})\label{eq:6Bis}\end{equation}
(the $+2$ appears as follows: the $\mathbf{B}_{0}$ for which our
argument doesn't work as is, and the very last $\mathbf{B}_{i}$ where
the interval $I_{i}$ might be cut off and we don't want to mess with
this problem).

We still need one calculation to overcome the condition $R(n+ik)_{1}>a_{i}$
in (\ref{eq:defBi}). Applying the same coupling argument as in the
beginning of the lemma shows that \[
\mathbb{P}(R(n+(i-1)k)_{1}<a_{i}-2n^{5/8})\leq Ce^{-cn^{1/8}}\ll C\exp (-c\sqrt{\log n}).\]
(\ref{eq:RIilarge}) gives, for $n$ sufficiently large as to satisfy
$2n^{5/8}<c\epsilon k^{3/4}$ \begin{equation}
\mathbb{P}(R(0,n+ik)_{1}>a_{i})>1-C\exp (-c\sqrt{\log n})\quad \forall i\label{eq:Rnikaionei}\end{equation}
and (\ref{eq:Rnikaionei}) holds for any $n$ if only $C$ is sufficiently
large as to make it trivial for smaller $n$'s. Summing over all $i$
we get \begin{equation}
\mathbb{P}(\exists i:R(n+ik)_{1}\leq a_{i})\leq \left\lceil \frac{n}{k}\right\rceil C\exp (-c\sqrt{\log n})\leq C\exp (-c\sqrt{\log n}).\label{eq:Rnikaialli}\end{equation}
Assuming that the events in (\ref{eq:Cirrelev}), (\ref{eq:6Bis})
and (\ref{eq:Rnikaialli}) did not happen we get that $R(2n)_{1}-R(n)_{1}\geq \alpha _{k}k\left(\left\lceil \frac{n}{k}\right\rceil -2\lambda -2\right)$
(note that the event of (\ref{eq:Rnikaialli}) implies, in particular,
that for the $2\lambda +2$ bad $i$'s, we still get that $R(n+k(i+1))_{1}-R(n+ik)_{1}>0$).
Define therefore \begin{equation}
\alpha _{n}=\alpha _{k}\left(1-\frac{2\lambda +2}{n/k}\right)\label{eq:defalphan}\end{equation}
and get (we now number some constants for clarity in the next part),\begin{eqnarray}
\lefteqn{\mathbb{P}(R(n)_{1}>0\textrm{ and }R(2n)_{1}<R(n)_{1}+\alpha _{n}n)\leq } &  & \label{eq:findCc1}\\
 & \qquad  & \leq C_{6}\exp (-c_{5}\sqrt{\log n})+C_{7}\exp (-c_{4}\lambda (1-o(1))\sqrt{\log n}).\nonumber 
\end{eqnarray}

We can now fix our parameters. Strangely enough, we start with $\lambda $,
and fix it so that $\left(\frac{11}{12}\right)^{\lambda -2}<1/2C_{7}$.
Next we fix $N_{1}(\epsilon )$ sufficiently large so that for all
$n>N_{1}(\epsilon )$ the $o(1)$ inside the second exponent is $<\frac{1}{\lambda }$,
and so that $\frac{n}{k}>4\lambda +4$, so that (\ref{eq:defalphan})
makes sense. Next we need to fix $c_{4}$. Let $N_{2}(x)$ be defined
by \[
\exp (-x\sqrt{\log N_{2}(x)})=\frac{11}{12}.\]
For $c_{4}$ sufficiently small, we would have \begin{equation}
C_{6}\exp (-c_{5}\sqrt{\log n})<\frac{1}{2}\exp (-c_{4}\sqrt{\log n})\quad \forall n>N_{2}(c_{4}).\label{eq:C2c2C1}\end{equation}
Further, for $c_{4}$ sufficiently small we would have $N_{2}(c_{4})>N_{1}$.
Fix $c_{4}$ sufficiently small so as to satisfy both conditions.
Hence, for all $n>N_{2}(c_{4})$ we get\begin{align}
C_{7}\exp (-c_{4}\lambda  & (1-o(1))\sqrt{\log n})\leq C_{7}\exp (-c_{4}(\lambda -1)\sqrt{\log n})\nonumber \\
 & \leq C_{7}\exp (-c_{4}(\lambda -2)\sqrt{\log N_{2}(c_{4})})\exp (-c_{4}\sqrt{\log n})\nonumber \\
 & \leq C_{7}\left({\textstyle \frac{11}{12}}\right)^{\lambda -2}\exp (-c_{4}\sqrt{\log n})<\frac{1}{2}\exp (-c_{4}\sqrt{\log n})\label{eq:C1c1c1}
\end{align}
and with (\ref{eq:C2c2C1}) and (\ref{eq:findCc1}) we get\[
\mathbb{P}(R(n)_{1}>0\textrm{ and }R(2n)_{1}<R(n)_{1}+\alpha _{n}n)\leq \exp (-c_{4}\sqrt{\log n})\quad \forall n>N_{2},\]
as required. Finally, for $n\leq N_{2}$, we only need to show that
\[
\mathbb{P}(R(2n)_{1}<R(n)_{1}+\alpha _{n}n)\leq \exp (-c_{4}\sqrt{\log n})\leq \frac{11}{12}.\]
However, setting $\alpha (n)=\frac{1}{2n}$ we get, from the coupling
of the excited random walk to the simple three-dimensional random
walk, as in the beginning of the lemma, that \[
\mathbb{P}(R(2n)_{1}\leq R(n)_{1})\leq \frac{5}{6}\]
(showing this for a three-dimensional simple random walk is a straightforward
calculation --- but if you really don't want to do it, replace $\frac{11}{12}$
everywhere with $1-c$ for some $c>0$) and we are done. 
\end{proof}
\begin{cor*}
Let $R$ be an $\epsilon $-excited walk of length $2n$ (starting
from $0$). Then \[
\mathbb{P}(R(2n)_{1}\leq c(\epsilon )n)\leq \exp (-c(\epsilon )\sqrt{\log n}).\]

\end{cor*}
\begin{proof}
Translate by $n+1$ so that $R$ starts from $(n+1,0,0)$. Use lemma
\ref{lem:excmain} with $\mathfrak{V}=\emptyset $ to get \[
\mathbb{P}(R(n)_{1}>0\textrm{ and }R(2n)_{1}<R(n)_{1}+\alpha _{n}(\epsilon )n)\leq C\exp (-c(\epsilon )\sqrt{\log n}).\]
Since $\alpha _{n}(\epsilon )\geq c(\epsilon )$, and since $R(n)_{1}>0$
always, we are done.
\end{proof}

\begin{proof}
[Proof of theorem \ref{thm:exc}]Use the last corollary for $n=2^{i}$
and get a sequence of event with probability $\leq \exp (-c\sqrt{i})$
hence by the Borel Cantelli lemma only a finite number of them occur,
and we are done.
\end{proof}

\end{document}